\documentclass{article}

\usepackage{arxiv}

\usepackage[utf8]{inputenc} 
\usepackage[T1]{fontenc}    
\usepackage{hyperref}       
\usepackage{url}            
\usepackage{booktabs}       
\usepackage{amsfonts}       
\usepackage{nicefrac}       
\usepackage{microtype}      
\usepackage{lipsum}
\usepackage{graphicx}
\usepackage{amsmath}
\graphicspath{ {./images/} }

\title{An interdisciplinary data-science approach to managing natural hazards risk}

\author{
 Cristobal Pais \\
  Department of Civil and Environmental Engineering\\  
  University of California Berkeley\\
  Berkeley, CA 94720 \\
  \texttt{cpaismz@berkeley.edu} \\
   \And
 Minho Kim \\
 Department of Landscape Architecture and Environmental Planning\\
 University of California Berkeley\\
 Berkeley CA 94720\\
  \texttt{mhk@berkeley.edu } \\
  \And
 Yanyan Xu\\
MoE Key Laboratory of Artificial Intelligence\\ 
AI Institute Shanghai Jiao Tong University\\
Shanghai 200240, China
  \texttt{yanyanxu@sjtu.edu.cn} \\
  \And
  John Radke \\
 Department of Landscape Architecture and Environmental Planning\\
 University of California Berkeley\\
 Berkeley CA 94720\\
  \texttt{ratt@berkeley.edu } \\
  \And
   Marta C. Gonzalez\\
  Department of Civil and Environmental Engineering\\  
  Department of City and Regional Planning\  
  University of California Berkeley\\
  Lawrence Berkeley National Laboratory\\ 
  Berkeley, CA 94720 \\
  \texttt{martag@berkeley.edu} \\
}

\begin{document}
\maketitle
\begin{abstract}
Natural hazard risk management is a demanding interdisciplinary task. It requires domain knowledge, integration of robust computational methods, and effective use of complex datasets. However, existing solutions tend to focus on specific aspects, data, or methods, limiting their impact and applicability. Here, we present a general data-driven framework to support risk assessment and policy making illustrating its usage in the context of fire hazard by integrating three unique datasets of fire behavior, street network, and census data for the whole state of California. We show that integrating spatial complexity by including a fire behavior layer and a socio-demographic layer changes the universal function observed in previous optimization frameworks that only work with the accessibility of facilities.
These results open avenues for the future development of flexible interdisciplinary frameworks in natural hazards management using complex large-scale data.
\end{abstract}

Natural hazards strength and impact in our society have been significantly exacerbated in the last decades\cite{iglesias2021risky,kreibich2014costing,walker2019increasing}. Attributed to climate change and human activities, the rise in temperature and perturbation of known climatic conditions have triggered dramatic events of unprecedented magnitude across the globe\cite{balch2017human,zou2021increasing}. Only in the last five years, we have witnessed some of the record breaking floods, wildfires, droughts, tornadoes, hurricanes, among others devastating natural hazards, causing irreparable human losses, damage to communities/ecosystems, and millionaire expenses in infrastructure\cite{floodingdv,mallakpour2015changing,fire2018top,spinoni2014world}.

As natural disasters are part of our life, they have been deeply studied by researchers and practitioners 
to (i) comprehensively understand them; (ii) predict their occurrence; and (iii) develop effective mitigation policies to proactively protect our communities\cite{mechler2016reviewing,rivera2006brief,rus2018resilience}. Several approaches involving complex methods from 
disciplines such as computer, environmental, and social sciences have been developed, in an attempt to cover these objectives\cite{seaberg2017review,zhou2018emergency}. With the advent of the big data and the data-science era, supported by computational advances, these approaches have been able to incorporate new elements\cite{yu2018big}. These include the capacity to process larger and more complex types of information to model the phenomena of interest and capture their uncertainty; develop new computationally demanding but more accurate methods with large-scale applicability for prediction and decision-making; and reduce the gap between theoretical and practical applications. 
Today, data are not a constraint. High-quality datasets, at a global scale, are within the reach of `one-click' from multiple platforms\cite{gorelick2017google,planetaryms}. Similarly, we can access a myriad of rich data resources from official entities \cite{cenpy,hifld,calfire,bennett2010openstreetmap,ubermovement,rollins2009landfire,sugarbaker20143d} (e.g., census, NASA, HIFLD, CALFIRE, EPA) to incorporate in the models/systems. This massive access to resources represents an opportunity but also a challenge: we need to identify the most relevant, consistent, and effective information to incorporate into the systems while keeping them as interpretable and useful as possible.  

Improvements in remote sensing technology, data accessibility/availability, and data-science methods have not been fully reflected in practice. For example, in California, the largest wildfires 
were recorded within the past seven years \cite{firestats,sun2019wildfire,los_2019} causing devastating damage in the affected regions with millions of losses. Meanwhile, nationwide fire protection costs have increased by 135\% since 1986\cite{usfiredept}. Similarly, the unusual flooding event registered in New York in 2021\cite{nytimesflood} caused by record atmospheric conditions, caught authorities and our communities unprepared, causing the death of 14 people in the area and a total of 43 people in the region. 
These, among other examples, warn us that successfully connecting the worlds of computational and environmental/spatial sciences into effective practical applications remains an elusive task.

Given the complexity of this challenge, this calls for the development of new interdisciplinary approaches able to: (i) facilitate the data collection and integration from any domain; (ii) exploit data to capture the dynamics of the phenomena with multiple methods; (iii) uncover relevant patterns; and (iv) design effective mitigation plans\cite{aghakouchak2018natural,hino2017managed}. All this needs to be supported by a flexible design to allow the continuous ingestion of updated information to obtain realistic and actionable solutions.
Advances in generating global-detailed environmental and social data and in their accessibility can be combined with state-of-the-art data-science methods and domain knowledge to unveil previously unseen patterns in environmental phenomena\cite{hampton2013big}. Exploiting large, detailed, and complex information -- previously unavailable -- with cutting-edge computational methods will be key to understanding and supporting the development of robust mitigation policies in the uncertain and challenging context of natural hazards risk management.




In an attempt to contribute in this direction, here we address three limitations of traditional risk assessment on natural hazards. First, we provide a data and domain-agnostic framework coupling any complex layers of information. Second, we integrate state-of-the-art interdisciplinary methods to provide a 360 $^{\circ}$ 
view of the system and managerial insights for policymakers. Finally, we apply robust spatial methods to integrate all components, focusing on the flexibility and scalability of the solution when applying large-scale computational methods.


To illustrate our framework, we study the risk exposure to fires in the state of California and analyze the effects of potential mitigation policies. We show the importance of incorporating multiple layers of information to capture the complexity of the system 
(see Methods). We identify the potential of each region to improve its current exposure under optimal resource allocation policies across the state. We also demonstrate the value of different intervention strategies and their trade-offs under multiple scenarios. Finally, we show how the effectiveness of the optimal policies varies across regions depending on their characteristics such as rural-urban distribution, street network coverage/accessibility, and fire potential. The incorporation of the environmental risk drastically changes the results of optimal allocation of facilities that are generally weighted by travel times and population only~\cite{xu2020deconstructing,gastner2006optimal,um2009scaling}

\section*{Multilayered risk assessment and policy generation}
The proposed framework brings together three distinct models to assess natural hazard risks and generate intervention policies. First, the natural hazard phenomenon (wildfires in our case) is simulated to characterize its behavior at multiple scales and capture its uncertainty and potential impact \cite{finney2006overview}. Then, we characterize the study area with multiple layers of information, relevant for the risk evaluation and mitigation strategies. Agnostic to input data, the framework allows the seamless integration of multiple sources and formats capturing sociodemographic, environmental, and topological (e.g., networks) data, among others, at different resolutions.  

\begin{figure*}
    \centering
    \includegraphics[scale=0.22]{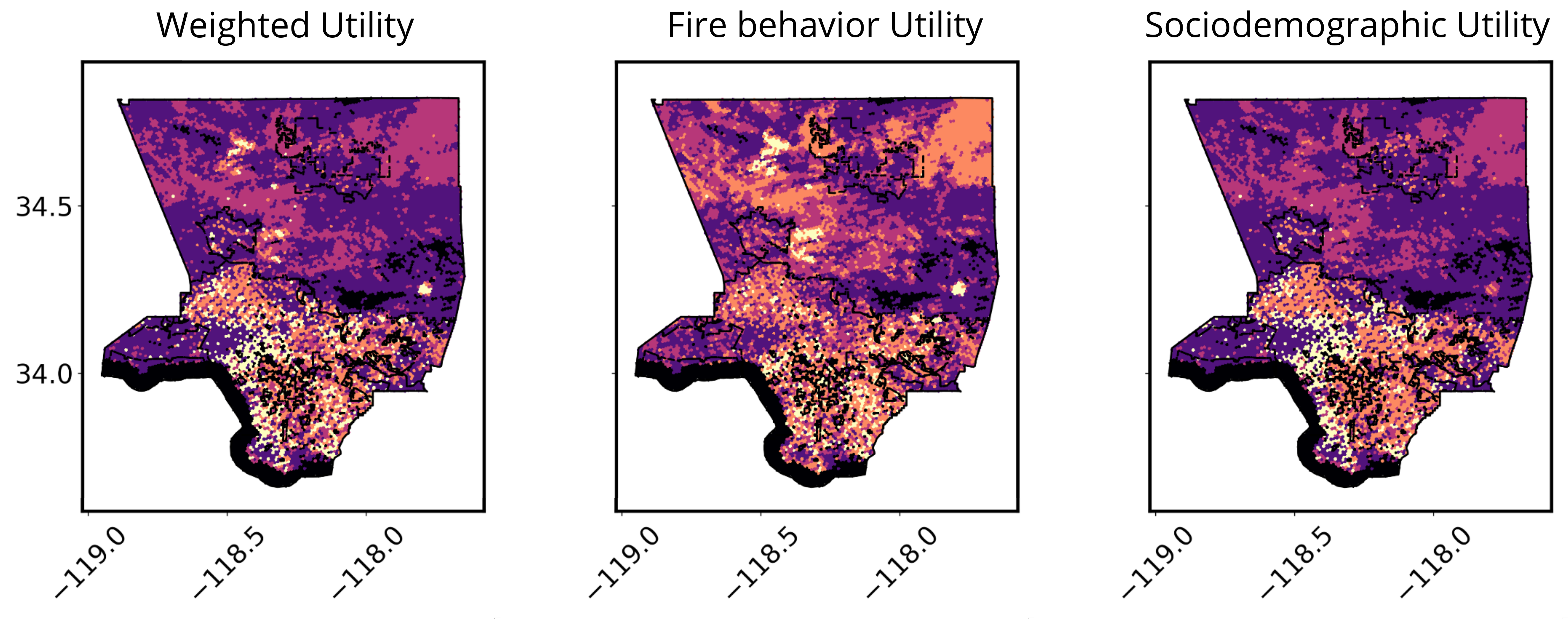}
    \caption{\textbf{Los Angeles county spatial distribution of combined features with different weights.} (left) Spatial distribution of the weighted average of fire behavior (FB) and sociodemographic (SD) features in the county. Areas with higher combined risk are represented by lighter colors. (center) Fire behavior dominated  (75\%) utility, highlighting those areas of the county with higher fire potential -- e.g., rural areas covered by grasslands located in the northern side of the county. (right) Sociodemographic dominated (75\%) utility highlighting densely populated and more expensive urban areas, mainly focused on the south of the county. }
    \label{fig:CALVaRDist}
\end{figure*}

Second, we consolidate data layers into a common resolution using an efficient spatial tessellation \cite{uberh3} obtaining a discrete map of the study area with aggregate statistics for each layer. We then apply a utility theory model \cite{elimbi2021quantifying} to generate scenarios using different weighted combinations of the layers, representing expectations/objectives of policy makers (Fig. \ref{fig:CALVaRDist}). Next, we construct a risk index (RI) highlighting the most vulnerable/exposed areas requiring potential interventions. This allows the evaluation of different scenarios, easily visualizing areas of agreement/disagreement among stakeholders while informing potential mitigation strategies under multiple objectives (see Methods).


Lastly, we use an optimization model to (i) prioritize areas of intervention according to the RI; (ii) determine the optimal allocation of resources/actions to improve the current scenario; and (iii) analyze the marginal contribution of adding extra resources in the areas of interest. Thus, we generate effective risk-mitigation policies according to a set of actions available to reduce the impact of the environmental hazard. Objectives, decisions, and constraints are modified accordingly (e.g., where to optimally locate fire stations to minimize wildfires exposure risk). Finally, we compare the initial state and the best-case solution, uncovering: improvement potential, optimal location of new/existing resources, and administrative level impact. For contrast, we compare the results of our framework (Fig. \ref{fig:PowerLaws}C) with the classical optimization that use travel times and population (Fig. \ref{fig:PowerLaws})B). We show that risk mitigation via facility allocation is slow and region dependent.



\section*{California and fire impacts}
California, home of more than 39 million people and the largest sub-national economy in the world with a \$3.4 trillion gross state product has experienced a notable increase in the annual burned area since 1972\cite{firestats}. This is driven primarily by longer and stronger fire seasons due to anthropogenic climate change effects \cite{williams2019observed}. Although fire suppression emergency funds have continuously increased in the last two decades -- 1.76 billion dollars in 2021 \cite{calfirefunds,usfiredept} -- the situation is far from ideal. Within the past seven years, it has recorded the largest wildfires in its history, causing devastating damage to communities and ecosystems, with losses valued in millions of dollars \cite{firestats,buechi2021long,li2021spatial}. Exacerbated by vast flammable landscapes (more than 50\% of the state) and a complex governance structure, days where people need to wear masks or use air filters to be able to breathe have become common every summer\cite{balmes2020changing}.

Los Angeles is the most populous (10M inhabitants) and the third densest county in California \cite{uscensusla}. Considered a crucial county in terms of valued assets\cite{calragelands} (observed median house value close to \$860,000 in 2020) and urban areas covering $\sim$27\% of the total county's surface, it has been severely impacted by wildfires in the last decade. With continuous urban growth and development, there is a concerning increase in (i) risk of ignition given the continuous growth in population and development in the wildland-urban interface (WUI)\cite{syphard2019factors}, and (ii) potential damage from catastrophic wildfires. Although having more than 400 fire stations, notable wildfires since 2009 have scorched large areas of communities, caused large evacuations, and destroyed hundreds to thousands of structures \cite{lacountyfirestats}. With more than 400,000 intersections, its street network is one of the most complex to navigate in the state, having approximately 1M different roads ($\sim$44\% of them residential).

Here, we simulate and analyze the impact of wildfires under extreme weather scenarios across the state to estimate worst-case expected fire behavior\cite{finney2006overview}. This is captured by the expected rate of spread (ROS) and fire intensity (FI) according to land vegetation and topography. From here, we evaluate and identify vulnerable areas by constructing a risk index combining sociodemographic (population, median house values, administrative/regional limits), environmental (land use, climatic, topography, and fire behavior), resources location (distribution of fire stations), and road network topology layers at a state-level. We consider a baseline case with the current resources allocation\cite{esri2007gis} (fire stations), accessibility conditions (road network), and fire behavior. We then explore the effectiveness of different mitigation policies -- relocation and addition of resources -- under multiple scenarios. 

\section*{Results}
\paragraph*{Risk-assessment and mitigation.}  

Focusing on a regional scale, we discuss the risk assessment and mitigation strategies using our framework in LA county. We analyze the distribution of different RIs, the expected improvement when relocating existing fire stations, and the trade-offs under this policy. 

From the simulations, we obtain a strong fire behavior characterized by an expected ROS of 2.4 km/h (max 5.5 km/h) and an expected FI of 0.2 KW/h (max 4.1 KW/h). This, is explained by the combination of favorable topographic, climatic, and vegetation conditions such as steep slopes covered by dry grass/shrubs,  becoming a threat and prone to ignition when temperatures are warm \cite{barbour2000north,keeley2002fire}.    

Observing the spatial distribution of three scenarios varying the importance of the sociodemographic and fire behavior features (Fig. \ref{fig:CALVaRDist}) we note that stronger fire behavior dominates the northern area of the county, mainly covered by shrubs and grasslands ($\sim$50\%) (middle panel). This is likely translated into faster, larger, and more intense fires than the ones we could expect in more urban areas, located in the south. 
Highly populated and dense cities such as Los Angeles, Santa Monica or Long Beach, with a urban-to-rural ratio near nine-to-one, can be found in these areas, with multiple and expensive structural assets exposed (right panel).

\begin{figure*}
    \centering
    \includegraphics[scale=0.13]{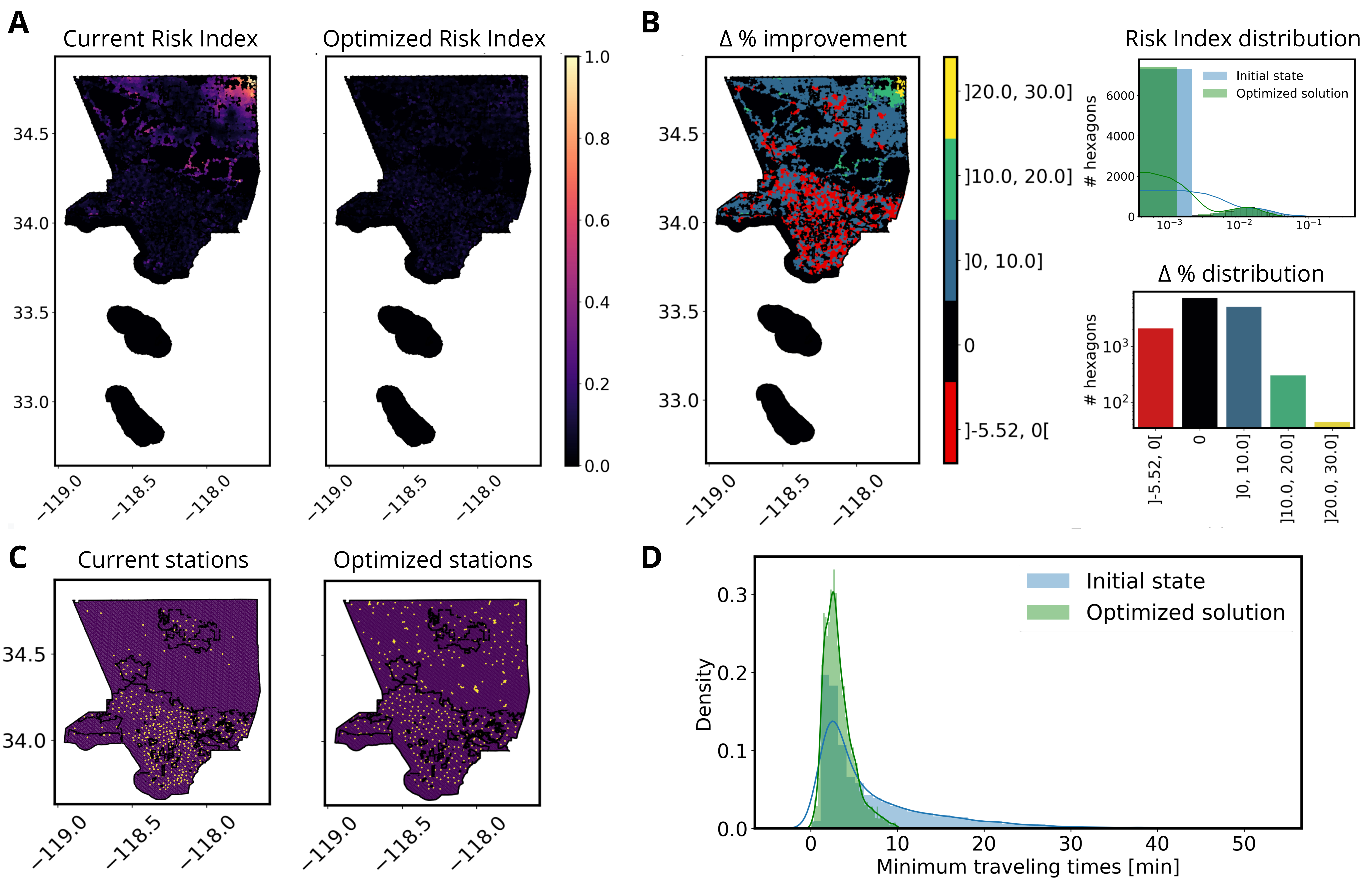}
    \caption{\textbf{Framework results in Los Angeles county.} (A) Spatial distribution of the current (left) and optimized (right) risk index across Los Angeles county. (B) Map of the county (left) indicating the risk index percentage change after optimizing the location of the existing fire stations inside the county. Higher values indicate better improvement potential with respect to the current state. Histograms of the current and optimized risk index (top-right) and the number of areas improved per bracket (bottom-right) across the county. (C) Comparison between the original (left) and the optimized (right) distribution of existing fire stations in the county, highlighted in yellow. (D) Distributions of the current (blue) and optimized (green) shortest traveling times in minutes between fire stations and each location of the county.}
    \label{fig:LACountyResults}
\end{figure*}

Introducing the location of fire stations in the street network to obtain RI immediately highlights clear accessibility issues to reach strong fire potential areas: fewer road options, slower travel times due to congestion, lower quality roads, and a reduced number of fire stations within 30 minutes (from $\sim$300 to only 10) compared to highly populated areas (see \textbf{Table S1} for detailed results). This leads to a long-tailed distribution of travel times in the northeastern region, with avg. and max values close to 19.6 and 45 minutes, respectively, almost 8.2 times the recommended response time by the National Fire Protection Association\cite{usfiredept}. We observe this in the distribution of RI (Fig. \ref{fig:LACountyResults}-A, left panel), noting a clear hot-spot at the north-east of the county, near Lancaster. Zooming into it, we find the dominance of flammable shrubs/scrubs (99.7\%) leading to above average RI values (max 0.33), mainly driven by faster fires (peaks of 4.6 km/h) while having below average population and exposed assets compared to the rest of the county. Therefore, the vulnerability of the area is driven by the mix of (i) an extremely dangerous fire potential associated with its vegetation/characteristics and (ii) a poor coverage and accessibility of existing fire stations due to lack of resources and the sparsity/complexity of the network to reach risky areas.

In contrast, the south is characterized by lower risk (max RI 0.27); slower fires (avg. 2 km/h); and average intensities.
Dominated by urban structures (55.4\%) providing home to a population close to 8.5M, shrubs (24.4\%), and water (12.8\%); the area experiences lower fire potential than its northern counterpart. However, on average, more expensive assets are exposed. Supported by great accessibility, it has multiple connections within the street network, concentrating almost 90\% of the existing stations. Therefore, average response times (3.4 minutes) are 82.6\% shorter than the northern region, reducing the overall vulnerability of the area.

After optimally relocating fire stations to minimize average RI across the county (Fig. \ref{fig:LACountyResults}-A), we see a significant decrease ($\sim$93\%) in the most vulnerable areas. This, triggered by an important redistribution of existing fire stations from most populated areas into rural and WUI northern locations (Fig. \ref{fig:LACountyResults}-C). This redistribution leads to an expansion of the average area covered by individual fire stations, increasing their population coverage by 7\%, on average. Interestingly, we note that stations may not be located right in the areas with higher risk. Better solutions could be found on nearby locations according to the street network: expected traveling times of the roads (correlated to their type and shape) and traffic directions are fundamental components to design the optimal policy.

Analyzing the trade-offs when improving the overall county situation, we see a global benefit in the RI distribution, reaching average and maximum values equal to 0.02 and 0.1, respectively -- a 50\% and 69\% improvement. For this, highly populated areas must accept worse RI values (Fig. \ref{fig:LACountyResults}-B, in red). Originally, these areas have access to a larger supply of fire stations, experiencing the shortest attention times in the county, thus leading to lower risk. However, this effect is overcome by the improvement of riskier areas: we observe a positive distribution between the current state and the optimal policy. From the panels, we note that RI is degraded for a total of 2,111 locations (avg. delta -0.08\%) while 5,455 locations improve their current situation (avg. delta 4\%, with peaks $\sim$30\%). We explain these results by the improvement in traveling times due to the stations relocation (Fig. \ref{fig:LACountyResults}-D): avg. and max traveling times are decreased by 56\% (from 7.3 to 3.2 minutes) and 80\% (from 50 to 10 minutes), respectively. Moreover, their distribution is compacted (i.e., smaller average, range, and maximum values), improving the county's risk exposure, on average.

\paragraph*{Objectives and policies sensitivity.}
We explore the impact of different objectives and constraints on optimal policies. First, we replace the objective to minimize (i) maximum RI values and (ii) a weighted function balancing average and maximum risk values (see Methods). Second, we optimize the location of new fire stations in the county, assuming that existing ones cannot be relocated. We further analyze their marginal contribution to minimizing risk and traveling times as we increase their number in the region. 

As expected, prioritizing the most vulnerable areas leads to aggressive re-localization of resources from urban into exposed rural zones, i.e., we obtain higher density of stations near the north-eastern hot-spot. This significantly reduces the maximum risk exposure and traveling times (by 84.3\% and 89.1\%, respectively). However, this aggressive redistribution has negative impacts, on average, as it prioritizes outliers over social optimum. Thus, we observe 2.3\% and 4.7\% larger RI and STTFS than the original policy. Focusing on the trade-offs, we note similar overall improvements across the county, but in this case, we reach higher values (max of 54\% instead of 30\%) since we prioritize the most vulnerable areas.

Simultaneously optimizing average and max risk across the county leads to obtaining negligible worse average results (less than 0.1\% worse) than the original policy. However, the new one decreases both the maximum RI and traveling times by 77\% and 80\%, respectively, about 26\% better than the initial model. To achieve this, 23 fire stations are relocated from the most vulnerable areas into strategic positions to improve their coverage according to the streets traffic and directions, and slightly relocating 47 stations within the densest urban areas to improve coverage near the WUI. As before, the new policy improves the regional situation. We note that this is the best among the three policies in terms of global reduction of average and maximum risk exposure, illustrating the importance of aligning model objectives and policy-maker goals. 

Fixing existing fire stations shifts the attention to where to optimally locate new facilities according to the policy-maker goals. In this case, we test minimizing the (i) mean and (ii) maximum risk values in the county. This approach allows us to identify of vulnerable locations, quantify potential improvements, prioritize efforts according to budget, and justify the investment in new resources for the community. 


Simulating the addition of 100 stations, we observe how they tend to be located near strong fire potential areas (north-eastern) (Fig. \textbf{S1}-(A,B)). Optimal locations are affected by roads' characteristics, transit directions, and accessibility of the area (e.g., the number of incoming streets). Thus, they may not be directly located in hot spots but within the neighborhood of the most vulnerable areas, reducing RI in the zone instead of just focusing on a few risky places. The latter is observed when optimizing max RI values. New stations are located right in the most vulnerable places. Thus, having an accurate street network representation is crucial: solutions may vary significantly depending on their composition. Second, we note an important decrease in the area covered and the total population served by each station of 37\% and 23\%, respectively. This leads to improved service levels (i.e., lower STTFS). Finally, new stations tend to not be located in dense urban areas as they are already well-supplied according to the policy goals.

Observing the impact of adding stations (Fig. \textbf{S1}-C), we note how initially they are located near high-vulnerable areas and thus, their marginal contribution to reducing RI is larger. For example, adding a unique station reduces max RI by 1.4\%. 
This effect is diluted as we increase the number of stations, reaching a saturation point where adding new resources becomes less effective, and thus, difficult to justify under strict budgets. In the simulation, we identify this threshold at $\sim$65 stations, where 
adding extra stations have a negligible contribution. At this threshold, we note a cumulative reduction of 46\% of max. RI, improving risk exposure for more than 60\% of the reachable locations. These analyses provide insightful support for policy-makers to (i) prioritize the use of resources; (ii) perform a robust cost-benefit analysis 
and (iii) quantitatively justify their decision-making process.

\paragraph*{Risk mitigation at state-level.}

When defining risk-mitigation policies, it is important to understand their impact and effectiveness at different aggregation levels. 
Having a holistic view, authorities could prioritize areas of interest according to (i) feasibility of implementing changes; (ii) improvement potential; and (iii) resources availability. We expand our analysis to the state of California. First, we analyze the distribution and potential risk improvement under optimal policies at two levels of aggregation: (i) cells and (ii) counties and core-based statistical areas (CBSA). Second, we discuss the trade-offs involved in the policy. Finally, we explore the characteristics of those areas with highest/lowest vulnerability values, identifying patterns and insights for policy-makers. 
\begin{figure*}
    \centering
    \includegraphics[scale=0.205]{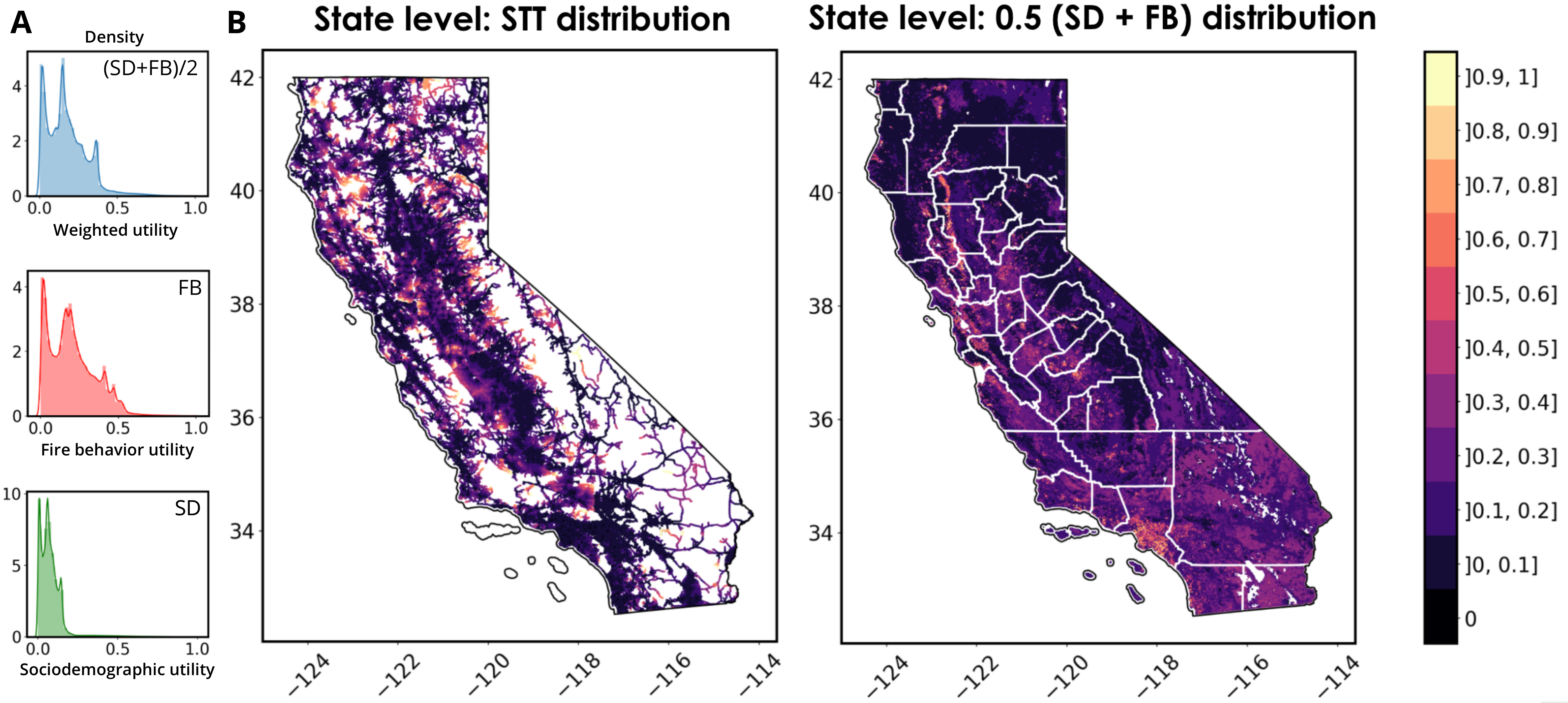}
    \caption{\textbf{Utility functions and traveling times spatial distribution in California.} (A) State-level density distribution plots for (i) the weighted utility function of the social planner (VaR utility); (ii) the fire behavior (FB) component of the utility function; and (iii) the sociodemographic (SD) component of the utility function. (B) State-level distribution of the shortest-traveling times between fire stations and each reachable point within the street network (left) and state-level spatial distribution of the normalized weighted utility including FB and SD components (right). Higher values are represented by lighter colors.}
    \label{fig:CALRIDist}
\end{figure*}

In Fig. \ref{fig:CALRIDist}-A, we note how the sociodemographic utility (green) tends to concentrate ($\sim$96\% of the state) below 0.25. This indicates that population and structure exposure risk is relevant only for a small proportion of the state: iconic densest urban areas like LA and the Bay area, representing only 4.4\% of cells. In contrast, fire behavior distribution (red) is wider, aligned with the presence of vast areas covered by flammable lands 
with peak speeds about 9 km/h (avg. 1.8 km/h). In panel B, we note areas experiencing accessibility difficulties, i.e., longer traveling times from stations. These areas suffer from: (i) lower street density; (ii) worse street conditions; and (iii) lower speed limits. These conditions lead to precarious situations: only one-third of the locations are reachable within the recommended response time\cite{usfiredept}. On average, it takes 14 minutes to reach any location from the nearest station and more than 10\% of them require +30 minutes, with peaks of 7.1 hours. At state-level (Fig. \ref{fig:CALRIDist}-B), we easily identify the most vulnerable locations combining high fire potential and assets exposure. For example, Los Angeles, aligned with our previous results: densely populated areas (south) and flammable areas (north) are the main drivers of vulnerability in the region. Joining both maps, we detect RI peaks of 0.85 in the most vulnerable areas and stations covering about 22,000 inhabitants on average.

\begin{figure*}
    \centering
    \includegraphics[scale=0.215]{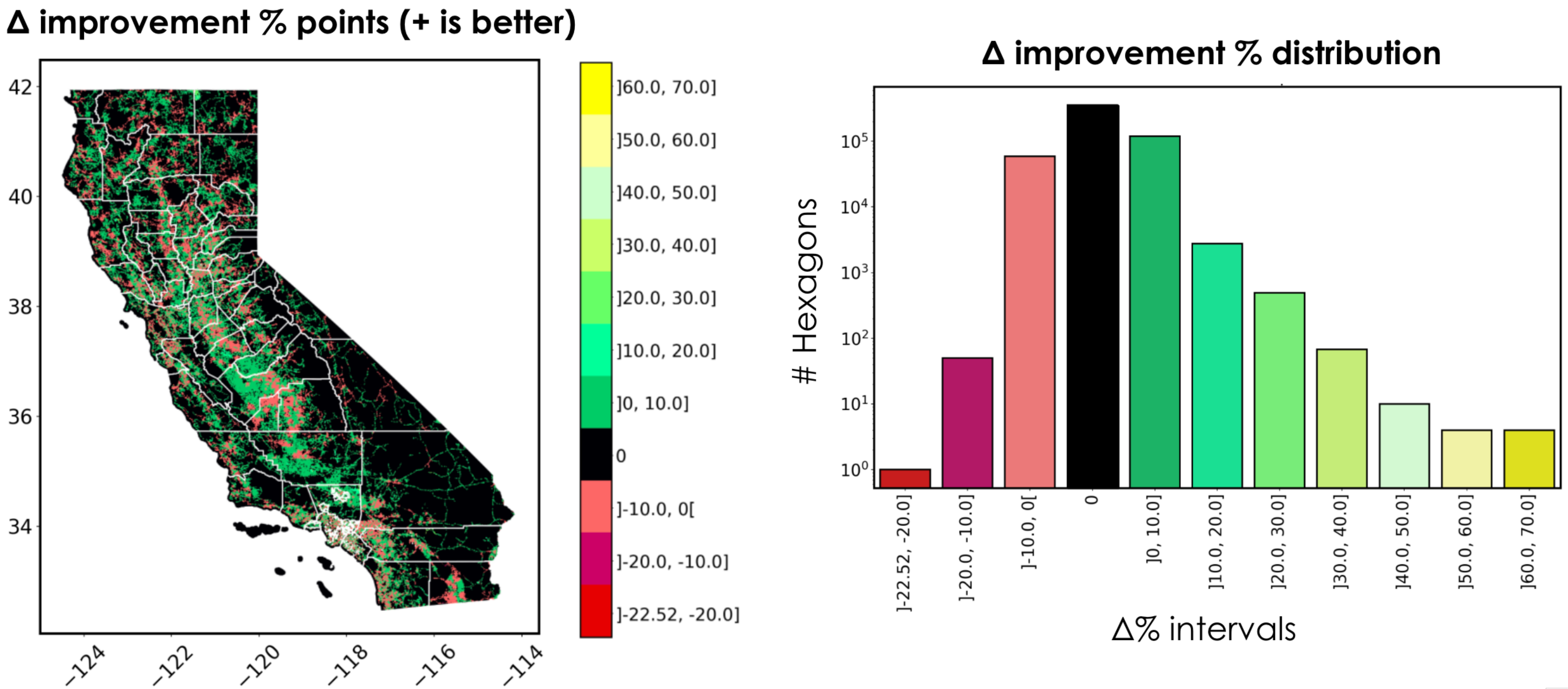}
    \caption{\textbf{State-level Risk Index improvement distribution at cell level.} Map of California divided into CBSAs and counties. (left) Risk index percentage change after optimizing the location of the fire stations across the whole state. Positive values indicate an improvement with respect to the current situation and negative values indicate a deterioration of the risk index in the area. (right) Distribution of the risk index percentage change indicating the total number of locations on each 10\% interval after optimizing the location of the fire stations to minimize the average risk index.}
    \label{fig:CALImprovement}
\end{figure*}

Results of the relocation policy (Fig. \ref{fig:CALImprovement}) indicate a global improvement towards risk-mitigation opportunities, representing 67\% of the road-accessible locations. Among this group, the majority has a potential risk improvement up to 40\% with most vulnerable zones reaching 70\%. 
Traveling times are improved by 77\% (avg. 3 minutes) and stations are better distributed, thus reducing supplied population by 37\%. This improves the overall risk exposure, reducing by 80\% the expected risk at state level.

Zooming out, we analyze the effectiveness of the policy across CBSAs and counties (Table \textbf{S1}). Today, the most exposed regions are Salinas, Orange County, and Vallejo (0.9 avg. RI). On the contrary, Amador County, Crescent City, and Susanville appear as the less vulnerable regions (0.13 avg. RI). To understand the main drivers of their risk, we analyze the contribution of the environmental, sociodemographic, and transportation layers on RI. We define three categories according to their relative contribution: (i) fire dominant; (ii) social dominant; and (iii) accessibility. We find El Centro, Merced, and Riverside as representatives of the first group. In these areas, the combination of highly-flammable land covers (shrubs/scrubs and crops); dangerous topographic characteristics boosting fires (up to 9 km/h); and limited accessibility; explains high average RI values ($\sim$0.6). In the second group, we find densely populated areas with relevant and expensive structural assets such as Orange and LA Counties, Sacramento-Roseville-Folsom, and San Francisco-Oakland (3.2M, 10M, 2.4M, and 5M inhabitants, respectively). These areas are characterized by dense, generally fast, and multi-directional congested street networks. The exposure of population/structures with flammable rural areas away from metropolis -- and thus, less accessible -- represent the main drivers of risk in these locations (RI of 0.85, 0.6, 0.5, and 0.3, respectively). An interesting case is San Francisco-Oakland: although dominated by densely urban locations and non-flammable areas (water), its transportation network suffers severe congestion events (e.g., near the bridges) triggering accessibility challenges. We find Trinity, Inyo, and Amador in the third group. These CBSAs experience the worst overall accessibility in the state. Their transportation networks are sparse, difficult/slow to navigate and suffer congestion. Thus, avg. and max traveling times are close to 46 minutes and 6 hours, respectively, leading to RI values close to 0.25.

\begin{figure*}
    \centering
    \includegraphics[scale=0.222]{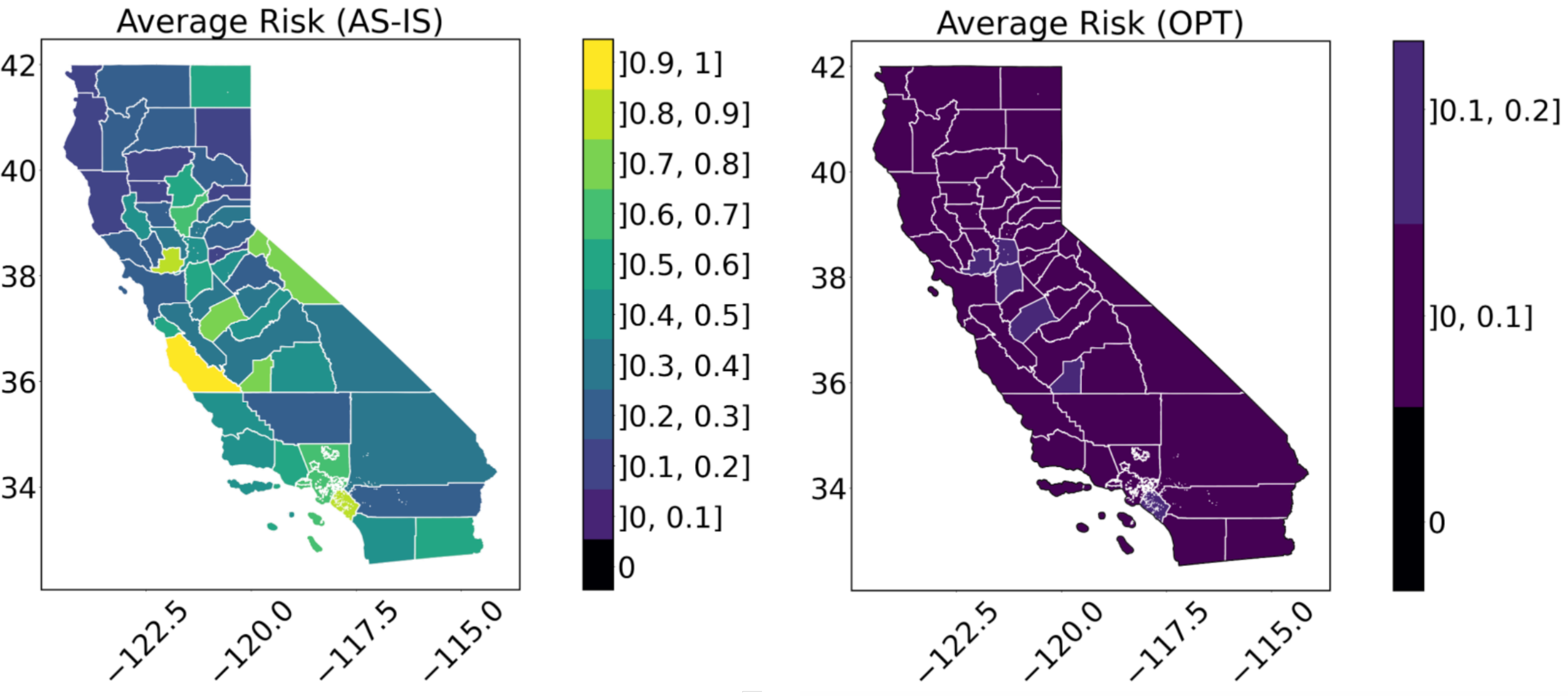}
    \caption{\textbf{State-level Risk index optimization results.} Map of California state divided into CBSA and counties. (left) Average risk index value distribution across California under the current conditions and fire stations location normalized by across all CBSA/counties. (right) Optimized risk index value distribution after allocating the existing fire locations in optimal locations.}
    \label{fig:CalAsISvsOpt}
\end{figure*}

Observing the normalized average RI values pre (AS-IS) and post (OPT) applying the policy (Fig. \ref{fig:CalAsISvsOpt}) indicates that the optimal policy leads to significant improvements: on average, RI and mean/max traveling times are reduced by 89\% and 72\%/51\%, respectively. Moreover, 48 regions can satisfy the recommended traveling time threshold under this policy, in contrast to only one in the current state (Orange county). This leads to improving average and peak state-level RI by 90\% and 84\%, respectively (Table \textbf{S1}). This pattern is explained by two factors. First, the average population covered by stations decreases due to a more effective assignment. Therefore, policy-makers can relocate resources to less accessible regions, commonly far from urban areas. Significant reductions in SFTTS to isolated areas are translated into lower risk exposure. Second, locations are selected by exploiting characteristics (e.g., transit directions) of the network to reach areas of interest faster. From this, Madera, Riverside-San Bernardino-Ontario, Inyo County, Alpine County, and Mono County are the areas with the greatest improvement potential in risk exposure (98\% improvement, on average).

Contrary, Sacramento-Roseville-Folsom, Glenn County, Stockton, Napa, and Merced experience the smallest improvement in RI (75\%) being Sacramento the most exposed region under the optimal policy. Similarly, Hanford-Corcoran, Orange County, and Yuba City experience the smallest improvement in traveling times ($\sim$45\%). This, is because these regions have: i) a lower stations/area ratio, leading to less flexible solutions; ii) sparse connectivity to peak traveling time areas; and iii) strong fire behavior. 

Although improving their overall risk exposure by 88\%, regions such as Stockton, Alpine County, Crescent City, and Calaveras County cannot significantly reduce peak traveling times ($\sim$4\% improvement) under the relocation policy. This highlights the need for resources in those areas. We observe similar patterns in El Centro, Santa Cruz-Watsonville, Colusa County, and Napa, still experiencing peak traveling times above one hour under the optimal policy.

\paragraph*{Facilities, optimal distributions, \& scaling laws.}
Previous works focusing on optimal facility distribution have studied the relationship between facility and population densities, identifying power functions to characterize AS-IS and optimal scenarios (for average distance) at the national scale \cite{um2009scaling,xu2020deconstructing}.  

From these studies, optimal power laws between the distribution of facilities and population show fitted exponents close to $2/3$ (see Methods). Moreover, it is shown that the optimal spatial distribution of facilities to minimize average traveling distance can be modeled using a general function derived from the analysis of multiple  cities\cite{xu2020deconstructing}. This allows policy-makers to plan for the optimal allocation of facilities without explicitly solving optimization models.

\begin{figure*}
    \centering
    \includegraphics[scale=0.095]{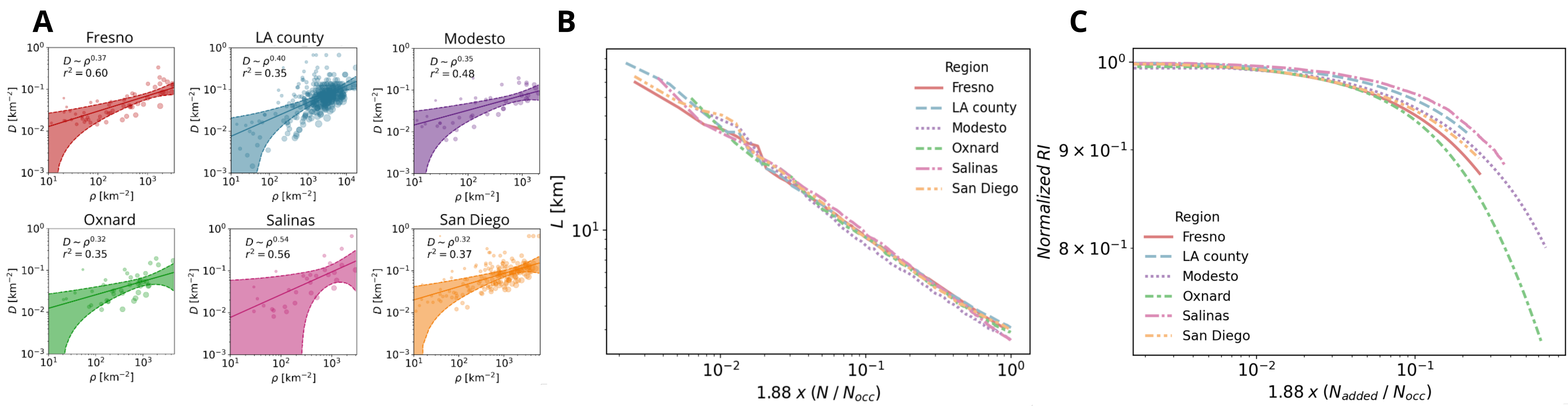}
    \caption{\textbf{Exploring the scaling laws of optimal distributions.} Sample regions analysis on facility distribution and impact of different objectives. (A) AS-IS power law analysis. For the actual distributions of facilities in several regions we observe the density of facilities $D$ scale with the population density as $D\sim\rho^{\beta}$ with $\beta \leq \beta_{opt} = 0.66$. (B) In the optimal scenario, optimizing purely for distance $L$, we see that the state’s fire facilities follow the previously reported\cite{xu2020deconstructing} function $L(\lambda N)$. (C) Our compound Risk index does not show a trivial behavior when adding more facilities in an optimal framework. This highlights the importance of incorporating environmental factors (i.e., fire behavior) in the optimization framework. A similar pattern is observed when optimizing the location of the existing facilities.}
    \label{fig:PowerLaws}
\end{figure*}

Inspired by these works, we study (i) AS-IS and optimal power laws; (ii) the validity of the general function when optimizing for average distance; and (iii) discuss the limitations of such a function when dealing with complex/compounded objectives such as our risk index, in multiple regions. Fig. \ref{fig:PowerLaws}-A shows how the AS-IS power laws indicate a theoretical non-optimal distribution of the facilities -- as expected -- across six sampled regions. This is reflected in the value of the $\beta$ exponents, different from the optimal $2/3$. Using our framework to optimize for a single objective as in previous works (i.e., traveling distance), we observe how: the optimal $\beta$ exponents tend towards the optimum value, and the validity of the proposed general function to model the optimal distribution of facilities. The latter is observed in Fig. \ref{fig:PowerLaws}-B, where curves represent the optimal avg. distance $L$, as a function of the scaled number of facilities by the total population, collapse into a single line for all regions. 

This universal shape is not observed when optimizing for a compound objective such as RI. Focusing on testing the validity of the general function when optimally allocating new facilities, we note (Fig. \ref{fig:PowerLaws}-C) how the relationship between optimal RI and the adjusted number of facilities is non-trivial: we cannot represent all curves with a unique function. We observe a similar pattern when solving the best-case scenario where existing fire stations are optimally relocated. These results highlight the important role played by the extra information layers (e.g., fire behavior in RI) incorporated into our framework: the presence of multiple objectives and interacting layers introduce a new level of complexity. This extra complexity in the decision-making process cannot be easily captured by approximation functions. It requires the explicit consolidation of all variables and data layers into an explicit optimization problem, as supported by the proposed framework. Moreover, the inclusion of extra features will likely lead to more complex results and patterns, making it extremely difficult for policymakers to use approximation functions without explicitly solving the problem under consideration.


\section*{Discussion}
The flexibility and interdisciplinary nature of the proposed framework allow for investigation questions incorporating all kinds of information while keeping the portability of the solutions. This allows policy-makers to freely characterize study areas with all their complexities, not limiting to specific features/formats. Coupling environmental, data-science, economic, and optimization methods, we provide a general approach to integrate unique knowledge from different areas in the challenging context of natural hazards. 

Policy-makers can use our framework as an evaluation and decision-making system to understand the impact of policies in different applications. For example, the impact of hazardous events such as flooding, hurricanes, or earthquakes could be incorporated via specific layers/simulations into the framework. Once all relevant layers are defined, alternative optimization models can be designed to incorporate domain knowledge when searching for optimal risk-mitigation policies. 

Illustrating its application in the context of wildfires in California, we (i) measure the vulnerability to fires according to a risk index capturing relevant characteristics of the area; (ii) evaluate mitigation policies under multiple scenarios; and (iii) uncover and analyze the trade-offs involved. Results reveal the opportunity to significantly improve wildfire risk exposure across the state.
We compare alternative mitigation policies, providing useful guidance about the impact of different actions, and quantitatively supporting policy-makers.  

Multiple layers of information can be easily added, allowing policymakers to evaluate their impact while enriching the policies. Layers such as the location of the electric distribution network (e.g., PG\&E), local traffic information, vehicle characteristics, presence of protected species or habitat niches, historical points of interest, and human mobility patterns could be used to evaluate and quantify their impact when generating complex risk-mitigation policies.

When it comes to the choice of a risk index, there is no standard solution. Multiple scenarios can be evaluated using different risk functions to characterize the study area, incorporate multiple objectives, and identify effective mitigation actions. Moreover, different regions may require alternative definitions to best represent their characteristics. Ultimately, this index should be chosen according to the objectives of stakeholders and policy-makers. Consolidating environmental, economic, and spatial constraints, drives the design of mitigation policies. Our framework supports the design, evaluation, and implementation of such a package of mitigation policies, unveiling relevant managerial insights by combining rigorous interdisciplinary methods into actionable policies.

\section*{Methods}
\label{mmethods}
The results presented in this paper are derived from a combination of frameworks and models using multiple data sources.

\paragraph*{Data sources.}
Sociodemographic data are obtained using the CenPy package, an interface used to automatically search and query the US Census API \cite{cenpy}. For this study, the data layers are queried from the American Community Survey in 2017 and  obtained per census tract (a total of 8,057 in California) for each county and core-based statistical areas (CBSA) as defined by the Census Bureau \cite{america2006us}. Data layers include median housing value (MHV), median household income (MHI), and total population (POP). 

Fire stations locations data are provided and consolidated from two main sources. (1) The CAL FIRE Facilities for Wildland Fire Protection 2020 dataset \cite{calfiredatasets}. It contains information about 1,255 state and local fire stations and other facilities -- owned and operated by CALFIRE\cite{calfire}, counties, or local government entities -- involved in wild land fire protection within the state of California. Besides the location of the station, records include facility name, facility type (based on funding source), funding, county, owner, staffing, address/phone, city, and zip code.  (2) Facilities location publicly available datasets from the Homeland Infrastructure Foundation-Level Data (HIFLD) \cite{hifld}. 

\paragraph*{Network topology and traveling times.}
The state-level road network is obtained from OpenStreetMap \cite{bennett2010openstreetmap} as polygonal files and explicit graph format using networkx package \cite{hagberg2008exploring}. It consists of 2,456,671 nodes (including intersections and dead ends) and 6,039,371 edges. When available, expected traveling speeds are inputted from publicly available data from a ride sharing company dataset\cite{ubermovement} and Google Maps to account for average traffic conditions. If no data is available (less than 5\% of the edges), legal traveling speed limits are used by default according to the type of road (e.g., highway or service road). 

Shortest traveling times (STT) paths between all the nodes of the network inside each county/CBSA (or relevant area of study) are calculated using the Dijsktra's algorithm \cite{dijkstra1959note} implementation in Python's igraph package \cite{csardi2006igraph}. In addition we calculate STT between the locations of the existing fire stations and all nodes of the street network (STTFS). From this, we also determine the current coverage and area of influence of each station, associating them to the corresponding shortest-traveling nodes. Isochrones are recorded for each node of the network for reference.

\paragraph*{Fire behavior modeling.}
A region of interest is modeled as a two-dimensional grid representing connectivity between hexagons/cells. Each cell represents a homogeneous area with similar characteristics (e.g., fuel type). Fire behavior data including rate of spread (ROS) and fireline intensity (FI) are generated using FlamMap \cite{finney2006overview}, a fire spread simulation model. The simulations use the Anderson's Fire Behavior model \cite{anderson1981aids} acquired from the Landscape Fire and Resource Management Planning Tools (LANDFIRE) \cite{rollins2009landfire} released in 2021. It consists of 13 fuel models grouped into four groups (grasslands, shrublands, timber, and slash) based on fuel type, fuel loading, size classes, moisture of extinction, and fuel bed depth. Terrain inputs including topography, elevation, slope, and aspect data are acquired from the USGS 3D Elevation Program \cite{sugarbaker20143d} released in 2017 at 1/3 arc second resolution. 

We obtain fire behavior variables across the study areas (counties and CBSA) using the 97th percentile fire weather and dead fuel moisture conditions (i.e., worst-case scenario). The model was designed for the severe period of the fire season and are generally useful for prediction fire behavior in such extreme scenarios \cite{anderson1981aids}. 

\paragraph*{Tessellation mapping.}
Data layers are integrated and consolidated into a two-dimensional mesh partitioned into hexagons of identical area. For this, we use the efficient Hexagonal Hierarchical Spatial Index (H3) \cite{uberh3}, which partitions the Earth into hexagons. This scheme has several advantages from spatial analysis and abstraction perspectives. First, the distance between adjacent hexagons' centers is always the same, simplifying and minimizing the error of spatial statistical analyses (in contrast to, e.g., a square raster system). Second, the discretization is extremely fast and allows dynamic resizing of the hexagons, enabling us to remap all features to different hexagons resolutions and see the impact of different aggregation levels in the generated policies. Moreover, different sections of the study area can be covered with different resolutions, making it efficient and effective when processing large areas where some sections do not require the same level of detail (e.g., areas with lower resolution data). Third, it seamlessly integrates with dynamic visualization frameworks (e.g., Kepler.gl\cite{keplergl}) where decision-makers can easily interact and simulate the impact of their decisions under multiple scenarios. 

In our experiments, we tested two resolutions levels characterized by an average hexagon area of 0.74 and 0.1 km$^2$; and an average hexagon edge length of 0.46 and 0.17 km,  respectively. These correspond to 539,133 and 3,773,919 hexagons covering the state of California (14,970 and 104,769 hexagons covering Los Angeles county, respectively). 

\textit{Point data.} Location data summarized in spatial points representing, e.g., fire stations and facilities data, are aggregated by summing their occurrence inside each hexagon obtaining frequency layers. A similar approach can be applied to event-type data (e.g., presence of a specific species of animal/vegetation the decision-maker is interested in).

\textit{Raster data.} Existing data formatted in traditional raster files (i.e., square grid) are re-sampled into the hexagons according to their data types. Land cover and fuel types data are summarized by their mode inside each hexagon while registering the proportion of each value in a separated layer for further analysis. Sociodemographic layers (POP, MHV, and MHI) 
are aggregated by summing and averaging their corresponding values. Other statistics at hexagonal level such as standard deviation, median, mode, and min/max values are also recorded for completeness in separated layers.  

\textit{Polygons data.} Vector type data representing historical fire scars and the street network spatial distribution in the area of study are intersected with the hexagons to estimate the proportion of the cumulative percentage of area burned from historical fires, the average number of nodes/edges of the street network, and the density of streets corresponding to each hexagon. For example, if a street intersects multiple hexagons, its length is proportionally divided among them. Similarly, if a historical fire scar covers multiple hexagons, each hexagon will only register the area that was contained inside its boundaries.   

\textit{Street networks data.} Each hexagon is associated with a unique node of the street network. This node is the closest to its centroid. Shortest traveling times between hexagons are represented by the STT between their associated street network nodes using the STT calculated as described in the Network topology and traveling times section. This approximation (proportional to the area of the hexagon) allows decision-makers to abstract the complexity of the street networks topology into a common tessellation including all relevant layers of information to generate realistic policies aligned with their expectations. We note that finer resolutions (i.e., smaller hexagons) will lead to more accurate STT between hexagons. However, from our experiments, no major differences (less than 1\%) where observed in the generated policies when comparing our results using hexagons of 0.74 km$^2$ and 0.1 km$^2$ resolution. This should be adjusted depending on the resolution, application, and level of detail required by the decision-maker.

Using this approach, we integrate all data layers into a unique mesh, making them comparable and abstracting the complexity of consolidating them for the decision-maker's analysis. 

\paragraph*{An integrating framework.}
We integrate multiple data layers to model the exposure/value of different sections of the landscape to wildfire risk. These data layers will be crucial to generate effective policies by capturing the reality of the study area. Based on the methodology presented in previous research \cite{elimbi2021quantifying}, we apply the following steps: (1) \textit{Scaling}. We scale all tessellated data layers into a common interval $[0,1]$ using an adequate transformation (e.g., logarithmic, linear). We note that these transformations are selected by decision-makers to represent the expected impact of the features on planning policies. As an example, areas with higher expected ROS are more vulnerable in case of a fire, therefore, an increasing function will capture this expectation. Moreover, we could define an upper bound for this layer, indicating that all areas with ROS values above a certain threshold will be mapped to the maximum value of the interval (i.e., equal to a unit), thus prioritizing them when defining the planning policy.  (2) \textit{Featurization}. Then, we combine scaled layers into meaningful features via convex combinations. This, to summarize them into relevant and interpretable categories containing correlated/complementary layers. For example, we could condense all fire related scaled layers into a new ``Fire behavior'' feature by assigning equal weights (i.e., all matter the same to the decision-maker). Similarly, population and mean house value can be condensed into a ``Demographic'' feature. Once we generate all features summarizing the main characteristics of the land, we define an (3) outcome variable. This outcome variable will quantify the relative importance of each section of the land according to the expectations of the decision-maker(s), captured by the transformations and weighted combinations applied in the previous steps. For example, a decision-maker could equally weight the ``Fire behavior'' and ``Demographic'' features to obtain an ``Exposure'' outcome variable.  

We note that, using different transformations in step (1) and weights in steps (2) and (3), decision-makers will be able to generate multiple scenarios and evaluate the impact on the planning policies. This methodology provides a flexible approach to i) easily integrate multiple layers of information; ii) evaluate the impact of including/excluding features; and iii) quantify the sensitivity of the policy when we modify their relative importance. 

In this work, five data layers (POP, MHV, ROS, FI, and STTFS) are condensed into a unique outcome variable that we name ``Risk index'' (RI). For this, we (1) apply increasing linear and non-linear scaling transformations to all layers, to represent a baseline scenario where higher values of the layers data lead to higher expected risks. Therefore, highly flammable places with poor service coverage (long traveling times) and high population density will be prioritized in the optimal policies (see Optimization models section). (2) Once scaled, two features condensing fire-related and sociodemographic features are generated. We define \textit{Fire behavior} (FB) and \textit{Sociodemographic} (SD) features as the equally weighted convex combinations of ROS and FI, and POP and MHV, respectively.   
\begin{eqnarray}
    FB &=& 0.5 \times (ROS + FI) \\
    SD &=& 0.5 \times (POP + MHV)
\end{eqnarray}

Then, (3) we define RI as the convex combination of FB and SD, multiplied by STTFS:
\begin{eqnarray}
    RI &=& 0.5 \times (FB + SD) \times STTFS 
\end{eqnarray}
From its definition, we note that RI $\in [0,1]$ will be higher in those places with i) strong fire potential (i.e., larger FB values); ii) more population and expensive assets exposed, using the MHV as a simple proxy to capture assets values; and iii) higher/slower service times from fire stations, i.e., less accessible areas. Therefore, higher values of RI will highlight more vulnerable areas given the described criteria - according to the selected scaling functions and weights - thus leading to obtain policies prioritizing the allocation of resources in those areas. Following its definition, we can quantify the impact of modifying the spatial location of our resources (fire stations) in the STTFS layer, directly impacting the distribution of RI. Therefore, in this work we can use the location of the fire stations as our main decision variables when generating optimal policies aiming to improve the distribution of RI. For example, depending on the decision-makers objectives, they could find the social optimal solution by minimizing the total value of RI in the area by relocating the existing fire stations, find the optimal number of stations to achieve specific target mean/max values of RI, or fully focus their attention on the riskiest areas and improve their situation by adding extra stations in the worst covered areas, among many other possibilities (see next section for details).

\paragraph*{Optimization models.}
We formulate an Integer (IP) and a Mixed Integer (MIP) Programming models to obtain planning policies. The first one allows decision makers to obtain the global optimal solution when all existing resources (fire stations) can be redistributed to optimize a specific objective function. This is an idealized scenario to highlight the best-case improvement policy used to quantify the improvement potential of an area and analyze trade-offs of relocating resources. In the second model, we evaluate the impact of optimally locating additional fire stations in the study area, assuming that existing ones are fixed. This model allows decision makers to allocate new resources effectively - in those areas where the impact is highest - as well as study the marginal utility (i.e., variation on the objective function) of adding a new fire station.


In this work, we illustrate the usage of our system with three objective functions: (1) minimize the average values of the outcome distribution; (2) minimize the maximum values of its distribution; and (3) a bi-objective model balancing the previous functions. To improve solving times, we use the optimal solution of (1) as a hot-start for (2) and the latter for (3), reducing average solving times by about 50\%. The formulation of the first model minimizing the average values of the outcome variable is the following:
\begin{eqnarray}
        ({IP_{RI}}) \hspace{.3cm} \min U &=& \sum_{i,j \in I} {RI}_i \times x_{(i,j)} \\ 
        s.t.\hspace{.5cm} \sum_{i \in I} x_{(i,j)} &=& 1 \hspace{1.0cm} \forall j \in I \hspace{2.3cm} \\
        x_{(i,j)}  &\leq&  y_{i} \hspace{0.9cm} \forall i, j \in I \hspace{2.5cm}\\
        \sum_{j \in I: j \neq i} x_{(i,j)}  &\geq&  y_{i} \hspace{0.9cm} \forall i \in I \hspace{2.5cm}\\
        \sum_{i \in I} y_i &=& S \hspace{1cm} \\
        x_{(i,j)}, y_{i} &\in&  \lbrace 0, 1 \rbrace \hspace{0.5cm} \forall i, j \in I 
\end{eqnarray} 
where $I :=$ set of potential locations (i.e., hexagons) connected to the street network; $S :=$ the total number of existing stations; $t_{(i,j)} :=$ the shortest traveling time between $i$ and $j$, incorporated in ${RI}_i$; ${RI}_i :=$ the value of the outcome variable of interest at hexagon $i$; and the decision variables: 
\begin{eqnarray} 
x_{(i,j)} &=& \text{1 if station in location \textit{i} \text{covers} \textit{j}} \\
y_{i} &=& \text{1 if a station is located in \textit{i}} 
\end{eqnarray}

In the formulation, Eq. (5) represent the set of covering constraints, i.e., every location is covered by a fire station. In Eqs. (6) and (7), we ensure that locations can only be covered by existing stations and that we open stations only in locations where other areas can be reached, respectively. Finally, we enforce the total number of stations available in the solution in Eq. (8). We note that explicit traveling time constraints associated with fire risk (e.g., arrive before $t$ minutes) can be easily added to the model to represent, e.g., maximum arrival times to points of interest such as hospitals/shelters, among other possibilities.

To model the second scenario (fixing existing stations), we define $E := \lbrace i \in I\; | \text{ station exists in i} \rbrace$ and replace $S$ in Eq. (8) by $\overline{S}:=  S + \delta$ where $\delta$ is a parameter representing the number of extra stations to be located/evaluated. This parameter is modified in our experiments to quantify the impact of adding new stations to the area of study, measuring their marginal contribution to the objective function and identify saturation points where adding extra resources is no longer useful. In addition, we fix the location of the existing stations by including the following constraint:
\begin{equation}
    y_i = 1 \hspace{1cm} \forall i \in E
\end{equation}

The second objective function (non-linear) is formulated by adding an auxiliary decision variable $\overline{RI}$ and a constraint to track the maximum outcome value of the solution to linearize the objective function:  
\begin{eqnarray}
(IP_{\overline{RI}}) \hspace{.3cm} \min U &=& \overline{RI}\\ 
s.t.\hspace{.5cm} && [5] - [9]\\
\overline{RI} &\geq& x_{(i,j)} \times t_{(i,j)} \times {RI}_j \hspace{0.25cm} \forall i, j \in I \\
\overline{RI} &\in&  \mathbb{R^+}
\end{eqnarray}
where $\overline{RI} :=$ maximum outcome value of the solution.

The third objective function (bi-objective) balances Eqs. [4] and [13] by weights $\alpha_1$ and $\alpha_2$ such that $\alpha_1 + \alpha_2 = 1$, selecting them according to the policy maker expectations.

 Models are implemented in Pyomo \cite{hart2017pyomo} modeling language in Python and solved with CPLEX v20.1 academic solver to optimality within a time limit of one hour using default settings for each study area.

\paragraph*{Sensitivity analysis.}

Our system provides the option to model the distribution of the data layers with different functions as well as weight them to construct a global metric based on the decision maker's most important considerations. We evaluate the impact of different weighting schemes at the policy generation level. We model three convex combinations of our different features for our experiments. First, we use an equally weighted combination of FB and SD. Then, we create two feature dominant RI outcome variables, where the selected dominant feature (i.e., FB or SD) contributes $75\%$ to the defined RI whereas the remaining layer is weighted by $25\%$. These three RI combinations (denoted RI, RIF, and RIS, respectively) are created and analyzed for each CBSA/county in our experiments. In addition, we experiment with different scaling functions on the five data layers, including exponential, logarithmic, and other non-linear functions. Interesting examples with significant impact on the optimal policies are discussed to highlight the flexibility of the system.

\paragraph*{Data availability.}
The datasets generated and/or analyzed during the current study are available from the corresponding author upon reasonable request. Several input datasets that support the findings of the study are available from different platforms and sources and restrictions may apply to the availability of such data. Processed datasets can be found at \url{http://www.github.com/humnet/C3_NaturalDisasters} as examples.

\paragraph*{Code availability.}
All code used to conduct this analysis is freely available at \url{http://www.github.com/humnet/C3_NaturalDisasters}.

\bibliographystyle{unsrt}  
\bibliography{references}  

\section*{Acknowledgments}
All authors acknowledge the support of C3.ai through the grant \textbf{Multiscale analysis for Improved Risk Assessment of Wildfires facilitated by Data and Computation} and the SERDP program of DoD through the grant \textbf{Networked Infrastructures under Compound Extremes (NICE)}.

\section*{Author contributions} C.P. conceived and designed the research; C.P. performed the experiments and the creation of the database; C.P. analyzed the data; M.K. organized the datasets; All authors contributed to the statistical analysis; C.P. and M.G. prepared the visualizations and data presentation; The first manuscript draft was written by C.P.; All authors discussed the results and contributed to the manuscript editing process; M.G. and J.R. provided language check and proofreading. M.G. supervised the research project.

\section*{Competing interests}
The authors declare no competing interests.

\section*{Additional Information}
\textbf{Supplementary Information} is available for this paper.
Correspondence and requests for materials should be addressed to Cristobal Pais.

\end{document}